\documentclass{amsart}
\usepackage[utf8]{inputenc}

\newcommand{\nirpdftitle}{Explicit Computations of Fundamental Classes}
\usepackage{import}
\newif\ifnirfiles
\nirfilesfalse
\newif\ifnirfancy
\makeatletter
\@ifclassloaded{amsart}{\nirfancyfalse}{\nirfancytrue}
\makeatother

\usepackage[margin=1in, marginparwidth=2cm]{geometry}
\usepackage[table,dvipsnames]{xcolor}
\usepackage{amsmath,amssymb,amsthm}
\usepackage{amsfonts}
\usepackage{asymptote}
\usepackage{cancel}
\usepackage{enumitem}
\usepackage{etoolbox}
\usepackage{graphicx}
\usepackage{mathdots}
\usepackage{mathtools} 
\usepackage{pgffor}
\usepackage{subfiles}
\usepackage{stmaryrd}
\usepackage{textcase}
\usepackage{tikz-cd}
\usepackage{todonotes}
\usepackage{xparse}
\usepackage{xr}

\ifnirfiles
	\usepackage{sty/fitch}
	\usepackage{sty/quiver}
\fi

\newtoggle{nirfancy}
\ifnirfancy
	\toggletrue{nirfancy}
\else
	\togglefalse{nirfancy}
\fi
\ifnirfancy
	\usepackage{footnotebackref}
	\usepackage{tocbibind}
	
	\usepackage{cabin}
	\usepackage[default]{cantarell}
\fi

\usepackage{fancyhdr}

\fancypagestyle{contentpage}{%
	\lhead{\textit{\rightmark}}
	\cfoot{\thepage}
}

\definecolor{wikipediadarkblue}{rgb}{0.023, 0.270, 0.676}
\providecommand{\nirpdftitle}{notes}
\ifnirfancy
	\hypersetup{
		colorlinks,
		citecolor=black,
		filecolor=black,
		linkcolor=wikipediadarkblue,
		urlcolor=wikipediadarkblue,
		pdftitle={\nirpdftitle},
		pdfauthor={Nir Elber}
	}
\else
	\AtBeginDocument{
		\hypersetup{
			citecolor=black,
			filecolor=black,
			linkcolor=wikipediadarkblue,
			urlcolor=wikipediadarkblue,
			pdftitle={\nirpdftitle},
			pdfauthor={Nir Elber}
		}
	}
\fi
\usepackage{hyperref}
\usepackage[capitalize]{cleveref}

\numberwithin{equation}{section}
\def\equationautorefname#1#2\null{%
	(#2\null)%
}

\newif\ifnirdebug
\nirdebugfalse

\ifnirdebug
	\usepackage{showkeys}
\fi

\usepackage{imakeidx}
\makeindex[intoc, title=List of Definitions]
\newcounter{indexcounterlabel}
\newcommand{\nirindexlabel}[1]{\label{indexentry:#1}}
\newcommand{\nirindex}[1]
{%
	\stepcounter{indexcounterlabel}%
	\nirindexlabel{\theindexcounterlabel}%
	\index{#1|hyperref[indexentry:\theindexcounterlabel]}%
}

\iftoggle{nirfancy}
{
	\usepackage{titlesec}
	\newif\iftoc

	\titleformat
		{\part} 
		[display] 
		{\cabin\bfseries\LARGE\scshape} 
		{\centering\LARGE Part \thepart} 
		{10mm} %
		{\centering\Huge} 
		[
			\thispagestyle{empty}
		] 
	\titlespacing*{\part}{0mm}{30mm}{30mm}
	\titleclass{\part}{top}

	\titleformat
		{\chapter} 
		[display] 
		{\cabin} 
		{} 
		{2in} %
		{
			\raggedleft
			\iftoc
				\vspace{2in}
			\else
				{\LARGE\textsc{Theme}~{\cantarell\thechapter}}\\ 
			\fi
			\Huge\scshape\bfseries
		} 
		[
			\vspace{-18pt}%
			\rule{\textwidth}{0.1pt}
			\vspace{0.0in}
		] 
	\titlespacing{\chapter}
		{0pt}
		{
			\iftoc
				-103pt+1in
			\else
				-127pt+1in
			\fi
		}
		{0pt}

	\titleformat
		{\section}
		{\Large\bfseries}
		{\thesection}
		{1em}
		{}
		[]
	\setcounter{tocdepth}{1}
}
{}

\usepackage{epigraph}

\usepackage{etoolbox}
\makeatletter
\newlength\epitextskip
\pretocmd{\@epitext}{\em}{}{}
\apptocmd{\@epitext}{\em}{}{}
\patchcmd{\epigraph}
	{\@epitext{#1}\\}
	{\vspace{-0.3in+20pt}\@epitext{#1}\\[\epitextskip]}
	{}
	{}
\makeatother

\renewcommand{\AA}{\mathbb A}
\newcommand{\RR}{\mathbb R}
\newcommand{\ZZ}{\mathbb Z}
\newcommand{\NN}{\mathbb N}
\newcommand{\QQ}{\mathbb Q}
\newcommand{\CC}{\mathbb C}
\newcommand{\FF}{\mathbb F}
\newcommand{\OO}{\mathcal O}

\newcommand{\floor}[1]{\left\lfloor{#1}\right\rfloor}

\newcommand{\op}[1]{\operatorname{#1}}
\newcommand{\mf}[1]{\mathfrak{#1}}

\newcommand{\id}{\operatorname{id}}

\newcommand{\into}{\hookrightarrow}
\newcommand{\onto}{\twoheadrightarrow}

\DeclareFontFamily{U}{dmjhira}{}
\DeclareFontShape{U}{dmjhira}{m}{n}{ <-> dmjhira }{}

\AtBeginDocument{%
	\mathchardef\ordinarycolon=\mathcode`:
	\mathcode`:="8000
}
\makeatletter
\newcommand{\coloncheck}{\@ifnextchar={\coloneqq\@gobble}{\ordinarycolon}}
\makeatother
\begingroup\lccode`~=`: \lowercase{\endgroup\let~}\coloncheck

\newcommand{\nirideasymbol}{%
	\begin{tikzpicture}[baseline=(x.base)]
		\draw[rounded corners=.01em] (-.05em,-1.07em)rectangle(.05em,.78em);
		\draw[fill=white,rounded corners=1.3] (0,.75em)--(.75em,0)--(0,-.75em)--(-.75em,0)--cycle;
		\draw[line width=0.2mm, line cap=round](-.4em,-1.07em)--(.4em,-1.07em);
		\node(x) at (0,0em) {};
		\node at (0,0em) {{\cabin\textbf{!}}};
	\end{tikzpicture}%
}
\newcommand{\nirwarnsymbol}{%
	\begin{tikzpicture}[baseline=(x.base)]
		\draw[rounded corners=.01em] (-.05em,-1.07em)rectangle(.05em,.78em);
		\draw[fill=white,rounded corners=1.3] (0,.75em)--(.75em,0)--(0,-.75em)--(-.75em,0)--cycle;
		\draw[line width=0.2mm, line cap=round](-.4em,-1.07em)--(.4em,-1.07em);
		\node(x) at (0,0em) {};
		\draw[
			line cap=but,
			line join=round,
			x=.5em,
			line width=0.5mm,
			y=1*(height("Z")-\pgflinewidth)*(1-sin(10)),
			rotate=-10,
			rounded corners=1.5pt,
		](-0.57, 0.57) -- (0.57, 0.57) -- (-0.57, -0.57) -- (0.57, -0.57);
	\end{tikzpicture}%
}

\usepackage{marginnote}
\makeatletter
\long\def\@mn@@@marginnote[#1]#2[#3]{%
	\begingroup
		\ifmmode\mn@strut\let\@tempa\mn@vadjust\else
			\if@inlabel\leavevmode\fi
			\ifhmode\mn@strut\let\@tempa\mn@vadjust\else\let\@tempa\mn@vlap\fi
		\fi
		\@tempa{%
			\vbox to\z@{%
				\vss
				\@mn@margintest
				\if@reversemargin\if@tempswa
						\@tempswafalse
					\else
						\@tempswatrue
				\fi\fi

					\llap{%
						\vbox to\z@{\kern\marginnotevadjust\kern #3
							\vbox to\z@{%
								\hsize\marginparwidth
								\linewidth\hsize
								\kern-\parskip
								\marginfont\raggedleftmarginnote\strut\hspace{\z@}%
								\ignorespaces#1\endgraf
								\vss
							}%
							\vss
						}%
						\if@mn@verbose
							\PackageInfo{marginnote}{xpos seems to be \@mn@currxpos}%
						\fi
						\begingroup
							\ifx\@mn@currxpos\relax\else\ifx\@mn@currpos\@empty\else
									\kern\@mn@currxpos
							\fi\fi
							\ifx\@mn@currpage\relax
								\let\@mn@currpage\@ne
							\fi
							\if@twoside\ifodd\@mn@currpage\relax
									\kern-\oddsidemargin
								\else
									\kern-\evensidemargin
								\fi
							\else
								\kern-\oddsidemargin
							\fi
							\kern-1in
						\endgroup
						\kern\marginparsep
					}%
			}%
		}%
	\endgroup
}
\makeatother
\newlist{listalph}{enumerate}{1}
\setlist[listalph,1]{label=(\alph*)}
\newlist{listroman}{enumerate}{1}
\setlist[listroman,1]{label=(\roman*)}

\usepackage{courier}
\usepackage{listings}
\usepackage{thmtools,thm-restate}

\usepackage[framemethod=TikZ]{mdframed}
\usepackage[framemethod=TikZ]{mdframed}
\usepackage{xpatch}
\makeatletter
\xpatchcmd{\endmdframed}
	{\aftergroup\endmdf@trivlist\color@endgroup}
	{\endmdf@trivlist\color@endgroup\@doendpe}
	{}{}
\makeatother

\iftoggle{nirfancy}
{
	\definecolor{nirlightblue}{HTML}{f7f7ff}
	\definecolor{nirdarkblue}{HTML}{1d1dbf}
	\declaretheoremstyle[
		mdframed={
			backgroundcolor=nirlightblue,
			linecolor=nirdarkblue,
			rightline=false,
			topline=false,
			bottomline=false,
			linewidth=2pt,
			innertopmargin=5pt,
			innerbottommargin=8pt,
			innerleftmargin=8pt,
			leftmargin=-2pt,
			skipbelow=2pt,
			nobreak
		},
		headfont=\normalfont\bfseries\color{nirdarkblue}
	]{nirbluebox}
}
{
	\declaretheoremstyle[bodyfont=\itshape]{nirbluebox}
}
\ifx\thechapter\undefined
	\declaretheorem[style=nirbluebox,name=Theorem,within=section]{thm}
\else
	\declaretheorem[style=nirbluebox,name=Theorem,within=chapter]{thm}
\fi
\declaretheorem[style=nirbluebox,name=Theorem,numbered=no]{thm*}
\declaretheorem[style=nirbluebox,name=Theorem,sibling=thm]{theorem}
\declaretheorem[style=nirbluebox,name=Theorem,numbered=no]{theorem*}

\declaretheorem[style=nirbluebox,name=Proposition,numbered=no]{prop*}
\declaretheorem[style=nirbluebox,name=Proposition,sibling=thm]{proposition}
\declaretheorem[style=nirbluebox,name=Proposition,numbered=no]{proposition*}

\declaretheorem[style=nirbluebox,name=Lemma,numbered=no]{lem*}
\declaretheorem[style=nirbluebox,name=Lemma,sibling=thm]{lemma}
\declaretheorem[style=nirbluebox,name=Lemma,numbered=no]{lemma*}

\declaretheorem[style=nirbluebox,name=Corollary,numbered=no]{cor*}
\declaretheorem[style=nirbluebox,name=Corollary,sibling=thm]{corollary}
\declaretheorem[style=nirbluebox,name=Corollary,numbered=no]{corollary*}

\iftoggle{nirfancy}
{
	\definecolor{nirlightred}{RGB}{250, 220, 220}
	\definecolor{nirdarkred}{HTML}{f40000}
	\declaretheoremstyle[
		mdframed={
			backgroundcolor=nirlightred,
			linecolor=nirdarkred,
			rightline=false,
			topline=false,
			bottomline=false,
			linewidth=2pt,
			innertopmargin=5pt,
			innerbottommargin=8pt,
			innerleftmargin=8pt,
			leftmargin=-2pt,
			skipbelow=2pt,
			nobreak
		},
		headfont=\normalfont\bfseries\color{nirdarkred}
	]{nirredbox}
}
{
	\declaretheoremstyle[bodyfont=\itshape]{nirredbox}
}

\declaretheorem[style=nirredbox,name=Conjecture,numbered=no]{conj*}

\declaretheorem[style=nirredbox,name=Question,numbered=no]{ques*}

\declaretheorem[style=nirredbox,name=Convention,numbered=no]{conv*}

\declaretheorem[style=nirredbox,name=Convention,numbered=no]{convention*}

\makeatletter
\declaretheorem[
	style=nirredbox,
	name=Idea,
	sibling=thm,
	postheadhook={\phantom{}\marginnote{\nirideasymbol}[-3pt]%
	\ifthmt@thisistheone
		\hspace{-2.2pt}%
	\else%
		{}%
	\fi}
]{idea}
\declaretheorem[
	style=nirredbox,
	name=Warning,
	sibling=thm,
	postheadhook={\phantom{}\marginnote{\nirwarnsymbol}[-3pt]%
	\ifthmt@thisistheone
		\hspace{-2.2pt}%
	\fi}
]{warn}
\makeatother

\iftoggle{nirfancy}
{
	\definecolor{nirdarkgreen}{HTML}{20660a}
	\definecolor{nirlightgreen}{HTML}{ebffeb}
	\declaretheoremstyle[
		mdframed={
			backgroundcolor=nirlightgreen,
			linecolor=nirdarkgreen,
			rightline=false,
			topline=false,
			bottomline=false,
			linewidth=2pt,
			innertopmargin=5pt,
			innerbottommargin=8pt,
			innerleftmargin=8pt,
			leftmargin=-2pt,
			skipbelow=2pt,
			nobreak
		},
		headfont=\normalfont\bfseries\color{nirdarkgreen}
	]{nirgreenbox}
}
{
	\declaretheoremstyle[bodyfont=\itshape]{nirgreenbox}
}

\declaretheorem[style=nirgreenbox,name=Axiom,numbered=no]{ax*}

\declaretheorem[style=nirgreenbox,name=Axiom,numbered=no]{axiom*}

\declaretheorem[style=nirgreenbox,name=Inventory,numbered=no]{inv*}

\declaretheorem[style=nirgreenbox,name=Inventory,numbered=no]{inventory*}

\declaretheorem[style=nirgreenbox,name=Notation,numbered=no]{notation*}

\declaretheorem[style=nirgreenbox,name=Definition,numbered=no]{defihelper*}

\iftoggle{nirfancy}
{
	\definecolor{nirdarkcyan}{HTML}{25805e}
	\definecolor{nirlightcyan}{HTML}{edfffb}
	\declaretheoremstyle[
		mdframed={
			backgroundcolor=nirlightcyan,
			linecolor=nirdarkcyan,
			rightline=false,
			topline=false,
			bottomline=false,
			linewidth=2pt,
			innertopmargin=5pt,
			innerbottommargin=8pt,
			innerleftmargin=8pt,
			leftmargin=-2pt,
			skipbelow=2pt,
			nobreak
		},
		headfont=\normalfont\bfseries\color{nirdarkcyan}
	]{nircyanbox}
}
{
	\declaretheoremstyle[bodyfont=\itshape]{nircyanbox}
}

\declaretheorem[style=nircyanbox,name=Example,numbered=no]{ex*}
\declaretheorem[style=nircyanbox,name=Example,sibling=thm]{example}
\declaretheorem[style=nircyanbox,name=Example,numbered=no]{example*}

\declaretheorem[style=nircyanbox,name=Non-Example,numbered=no]{nex*}

\declaretheorem[style=nircyanbox,name=Non-Definition,numbered=no]{ndefi*}

\declaretheorem[style=nircyanbox,name=Exercise,numbered=no]{exercise*}

\declaretheorem[style=nircyanbox,name=Exercise,numbered=no]{exe*}

\iftoggle{nirfancy}
{
	\definecolor{nirlightbrown}{RGB}{252,240,235}
	\definecolor{nirdarkbrown}{HTML}{7a5342}
	\declaretheoremstyle[
		mdframed={
			backgroundcolor=nirlightbrown,
			linecolor=nirdarkbrown,
			rightline=false,
			topline=false,
			bottomline=false,
			linewidth=2pt,
			innertopmargin=5pt,
			innerbottommargin=8pt,
			innerleftmargin=8pt,
			leftmargin=-2pt,
			skipbelow=2pt,
			nobreak
		},
		headfont=\normalfont\bfseries\color{nirdarkbrown}
	]{nirbrownbox}
}
{
	\declaretheoremstyle[bodyfont=\itshape]{nirbrownbox}
}
\declaretheorem[style=nirbrownbox,name=Remark,sibling=thm]{remark}
\declaretheorem[style=nirbrownbox,name=Remark,numbered=no]{remark*}

\declaretheorem[style=nirbrownbox,name=Quote,numbered=no]{quot*}

\makeatletter
\def\thmt@setheadstyle#1{%
	\thmt@style@headstyle{%
		\def\NAME{\the\thm@headfont ##1}%
		\def\NUMBER{\bgroup\@upn{##2}\egroup}%
		\def\NOTE{\if=##3=\else\bgroup\thmt@space\the\thm@notefont(##3)\egroup\fi}%
		\def\NIRNOTE{\if=##3=\else##3\fi}%
	}%
	\def\thmt@tmp{#1}%
	\@onelevel@sanitize\thmt@tmp
	\ifcsname thmt@headstyle@\thmt@tmp\endcsname
		\thmt@style@headstyle\@xa{%
			\the\thmt@style@headstyle
			\csname thmt@headstyle@#1\endcsname
		}%
	\else
		\thmt@style@headstyle\@xa{%
			\the\thmt@style@headstyle
			#1%
		}%
	\fi
}
\renewenvironment{thmt@restatable}[3][]{%
	\thmt@toks{}
	\stepcounter{thmt@dummyctr}
	\long\def\thmrst@store##1{%
		\@xa\gdef\csname #3\endcsname{%
			\@ifstar{%
				\thmt@thisistheonefalse\csname thmt@stored@#3\endcsname
			}{%
				\thmt@thisistheonetrue\csname thmt@stored@#3\endcsname
			}%
		}%
		\@xa\long\@xa\gdef\csname thmt@stored@#3\@xa\endcsname\@xa{%
			\begingroup
			\ifthmt@thisistheone
				\thmt@rst@storecounters{#3}%
			\else
				\@xa\protected@edef\csname the#2\endcsname{%
					\thmt@trivialref{thmt@@#3}{??}}%
				\ifcsname r@thmt@@#3\endcsname\else
					\G@refundefinedtrue
				\fi
				\@xa\let\csname c@#2\endcsname=\c@thmt@dummyctr
				\@xa\let\csname theH#2\endcsname=\theHthmt@dummyctr
				\let\label=\thmt@gobble@label
				\let\index=\@gobble
				\let\ltx@label=\@gobble
				\def\thmt@restorecounters{}%
				\@for\thmt@ctr:=\thmt@innercounters\do{%
					\protected@edef\thmt@restorecounters{%
						\thmt@restorecounters
						\protect\setcounter{\thmt@ctr}{\arabic{\thmt@ctr}}%
					}%
				}%
				\thmt@trivialref{thmt@@#3@data}{}%
			\fi
			\ifthmt@restatethis
				\thmt@restatethisfalse
			\else
				\csname #2\@xa\endcsname\ifx\@nx#1\@nx\else[{#1}]\fi
			\fi
			\ifthmt@thisistheone
				\label{thmt@@#3}%
			\fi
			##1%
			\csname end#2\endcsname
			\ifthmt@thisistheone\else\thmt@restorecounters\fi
			\endgroup
		}
		\csname #3\@xa\endcsname\ifthmt@thisistheone\else*\fi
		\@xa\end\@xa{\@currenvir}
	}
	\thmt@collect@body\thmrst@store
}{%
}
\makeatother
\iftoggle{nirfancy}
{
	\declaretheoremstyle[
		mdframed={
			backgroundcolor=nirlightgreen,
			linecolor=nirdarkgreen,
			rightline=false,
			topline=false,
			bottomline=false,
			linewidth=2pt,
			innertopmargin=5pt,
			innerbottommargin=8pt,
			innerleftmargin=8pt,
			leftmargin=-2pt,
			skipbelow=2pt,
			nobreak
		},
		headfont=\normalfont\bfseries\color{nirdarkgreen},
		headformat={\NAME\ \NUMBER\NOTE{\nirindex{\NIRNOTE}}}
	]{nirgreenboxindexed}
}
{
	\declaretheoremstyle[bodyfont=\itshape]{nirgreenboxindexed}
}

\declaretheorem[style=nirgreenboxindexed,name=Definition,numbered=no]{defi*}
\declaretheorem[style=nirgreenboxindexed,name=Definition,sibling=thm]{definition}
\declaretheorem[style=nirgreenboxindexed,name=Definition,numbered=no]{definition*}
\numberwithin{equation}{section}
\usepackage[backend=biber,
    style=alphabetic
]{biblatex}
\title{Explicit Computations of Fundamental Classes}
\author{Nir Elber}
\thanks{Research partially supported by NSF DMS-1840234.}
\date{\today}
\usepackage{graphicx}
\begin{document}

\maketitle

\begin{abstract}
	\noindent We use the techniques of group cohomology to give explicit computations of the local fundamental class. As an application, we discuss how to compute the Tate canonical class for the extension $\QQ(\zeta_{p^\nu})/\QQ$, where $p^\nu$ is an odd prime power.
\end{abstract}

\setcounter{tocdepth}{4}
\tableofcontents

\section{Introduction} \label{sec:intro}

Given a Galois extension of local fields $L/K$, the local fundamental class is a central construction in algebraic number theory because it encodes the Artin reciprocity map
\[\theta_{L/K}\colon K^\times/\op N(L^\times)\to\op{Gal}(L/K)^{\mathrm{ab}}.\]
See, for example, \cite[Theorem~2.2]{local-cft}. The goal of this paper is to provide an explicit computation of the local fundamental class of a Galois extension of local fields $L/K$ when $\op{Gal}(L/K)$.

The main idea for the computation is to use techniques motivated by abelian crossed products to appropriately abbreviate the data of $2$-cocycle. Our main result is as follows.
\begin{theorem} \label{thm:we-can-compute-fund-class}
    Fix a finite abelian extension of local fields $L/K$ with Galois group $G$. We show how to compute a $2$-encoding tuple representing the local fundamental class ($2$-encoding tuples will be defined in \Cref{def:encoding-tuple}). In particular, if $G$ is generated by $m$ elements, then we only require $m+m^2$ elements of $L^\times$ to describe the local fundamental class.
\end{theorem}
As an aside, at the cost of using more elements to describe the local fundamental class, one is able to use our techniques to remove the hypothesis that $L/K$ is abelian; see \Cref{rem:general-local-fund}.

The given presentation of the local fundamental class has a number of benefits. For example, it turns out that one can make the local fundamental class sufficiently explicit when $L/K$ has controlled ramification. Here is the more precise statement.
\begin{proposition} \label{prop:tame-tuple}
    Let $L/K$ be a finite Galois extension of tamely ramified extensions with ramification index $e$ and inertial degree $f$, where $K$ has residue field of order $q$. Further, suppose that $e\mid q-1$. To compute the fundamental class of $L/K$, we require the following data.
    \begin{itemize}
        \item Let $\sigma_0$ be the Frobenius automorphism of the maximal unramified subextension of $L/K$, extended to $\op{Gal}(L/K)$.
        \item Let $\pi$ be a uniformizer of $K$. Then we find a uniformizer $\gamma\in L^{\langle\sigma_0\rangle}$ such that $\gamma^e=\pi$ so that $L^{\langle\sigma_0\rangle}=K(\gamma)$.
        \item We find a primitive $\left(q^f-1\right)$st root of unity $\zeta_{q^f-1}$
        \item With $\zeta_e\coloneqq\zeta_{q^f-1}^{\left(q^f-1\right)/e}$, we define the automorphism $\sigma_1\in\op{Gal}(K(\gamma)/K)$ by $\sigma_1\colon\gamma\mapsto\zeta_e\gamma$. It turns out $\sigma_1$ generates this Galois group. Now, extend $\sigma_1$ to $L$ by fixing the unramified part.
    \end{itemize}
    Then if $e>1$, the fundamental class $u_{L/K}\in H^2(\op{Gal}(L/K),L^\times)$ is represented by the $2$-encoding tuple
    \[(\alpha_0,\alpha_1)=\left(\gamma,\zeta_{q^f-1}^{-1}\right)\qquad\text{and}\qquad
    \begin{pmatrix}\beta_{00} & \beta_{01} \\ \beta_{10} & \beta_{11}\end{pmatrix} = \begin{pmatrix}1 & \zeta_{q^f-1}^{-(q-1)/e} \\ \zeta_{q^f-1}^{(q-1)/e} & 1\end{pmatrix}.\]
    If $e=1$, then the same statement holds, ignoring the data at index $0$.
\end{proposition}
Another benefit to having access to a concrete local fundamental class is that it becomes clearer how the fundamental class encodes $\theta_{L/K}$, allowing us to more or less recover Dwork's \cite[Theorem~1]{explicit-local-recip-dwork} via purely cohomological techniques. Combining the above two observations allows us to compute $\theta_{\QQ_p(\zeta_p)/\QQ_p}$ without Lubin--Tate theory; see \Cref{ex:basic-ramified-artin-map}.

However, our presentation of the local fundamental class also makes it clear what data cannot be recovered from $\theta_{L/K}$; see \Cref{rem:cannot-get-betas-from-artin}. Though these data might appear without use, they are needed to define the Tate canonical class as in \cite{tate-torus}. Thus, as our last application is the following.
\begin{theorem} \label{thm:we-can-tate-class}
    We describe how to compute Tate canonical class of $\QQ(\zeta_{p^\nu})/\QQ$ when $p^\nu$ is an odd prime power.
\end{theorem}
There is already some literature in understanding the Tate canonical class (see, for example, \cite{debeerst-cft-algorithms} or \cite{chinburg-tate-class}), largely to understand the equivariant Tamagawa number conjectures. Notably, the main idea of our computation is similar to \cite[Section~4.5]{debeerst-cft-algorithms}.

\subsection{Overview}
In \cref{sec:local}, we use the theory of abelian crossed products to write down local fundamental classes of abelian extensions. In particular, \cref{subsec:group-cohom} recalls the cohomological tools we will need, culminating in an explicit computation of Restriction--Inflation in \Cref{prop:explicitresinf}. The computation of the local fundamental class is given in \cref{subsec:local-computation}, and the aforementioned applications follow it.

Lastly, in \cref{sec:global}, we compute the Tate canonical class of $\QQ(\zeta_{p^\nu})/\QQ$ when $p^\nu$ is an odd prime-power. Indeed, \cref{subsec:tate-class-defi} recalls the construction of the Tate canonical class, and the rest of the section is devoted to its computation.

\subsection{Acknowledgements}
This research was conducted at the University of Michigan REU during the summer of 2022. The author would especially like to thank his advisors Alexander Bertoloni Meli, Patrick Daniels, and Peter Dillery for their eternal patience and guidance, in addition to a number of helpful comments on earlier drafts of this article. Without their advice, this project would have been impossible. The author would also like to thank Maxwell Ye for a number of helpful conversations and consistent companionship. Without him, the author would have been left floating adrift and soulless.

\section{Local Fundamental Classes} \label{sec:local}
The goal of this section is to state \Cref{thm:we-can-compute-fund-class} precisely and provide its proof. We will use the theory of abelian crossed products, as presented in \cite{cohom-abelian-crossed}. Notably, our approach does not follow \cite[Theorem~XIII.2]{local-fields}; instead, our approach is purely cohomological.

\subsection{Group Cohomology} \label{subsec:group-cohom}
In this subsection, we pick up some cohomological tools that we will need for our computation.

\subsubsection{Abelian Crossed Products}

In this subsubsection, we recall the theory of abelian crossed products \cite{abelian-crossed}, recast in terms of group cohomology. Let $G$ be a finite abelian group, where $G=\bigoplus_{i=1}^m\langle\sigma_i\rangle$, where $\sigma_i\in G$ has order $n_i$. The idea here is to store data about $2$-cocycles in $Z^2(G,A)$ by only focusing on the generators $\sigma_i$. For brevity, we define the elements
\[T_i\coloneqq\sigma_i-1,\qquad N_i\coloneqq\sum_{k=0}^{n_i-1}\sigma_i^k,\qquad\text{and}\qquad\sigma^{(a)}\coloneqq\sum_{k=0}^{a-1}\sigma^k\]
for each $i$ and $\sigma\in G$ and $a\in\NN$.
\begin{definition} \label{def:encoding-tuple}
    Let $A$ be a $G$-module. Then a \textit{$2$-encoding tuple} is a pair $(\alpha,\beta)\in A^m\times A^{m\times m}$ of an $m$-tuple $\alpha$ and an $m^2$-tuple $\beta$ satisfying the following relations for each $i$ and $j$.
    \begin{itemize}
        \item $\beta_{ii}=0$.
        \item $\beta_{ij}=-\beta_{ji}$.
        \item $T_i\alpha_j=N_j\beta_{ij}$.
        \item $T_i\beta_{jk}+T_j\beta_{ki}+T_k\beta_{ij}=1$.
    \end{itemize}
\end{definition}
The notion of $2$-encoding tuples will be exactly what we need to discuss $2$-cocycles. An explanation of these relations can be found in \cite[Lemma~1.2]{abelian-crossed} and \cite[p.~423]{cohom-abelian-crossed}. The following proposition explains how to take a $2$-encoding tuple to a $2$-cocycle.
\begin{proposition}[{\cite[Theorem~2.1]{abelian-crossed-drozd}}] \label{prop:tuple-to-cocycle}
    Let $A$ be a $G$-module, and let $(\alpha,\beta)$ be a $2$-encoding tuple. Then the $2$-cochain
    \[(g,h)\mapsto\Bigg(\sum_{i=1}^m\Bigg(\prod_{k=1}^{i-1}\sigma_k^{a_k}\sigma_k^{b_k}\Bigg)\alpha_i^{\floor{(a_i+b_i)/n_i}}\Bigg)+\Bigg(\sum_{i>j}\Bigg(\prod_{k=1}^{i-1}\sigma_k^{a_k}\Bigg)\Bigg(\prod_{k=1}^{j-1}\sigma_k^{b_k}\Bigg)\sigma_i^{(a_i)}\sigma_j^{(b_j)}\beta_{ij}\Bigg)\]
    is a $2$-cocycle, where $g=\prod_{i=1}^m\sigma_i^{a_i}$ and $h=\prod_{i=1}^m\sigma_i^{b_i}$ such that $0\le a_i,b_i<n_i$ for each $i$. In fact, all cohomology classes are represented by a $2$-cocycle coming from some $2$-encoding tuple $(\alpha,\beta)$.
\end{proposition}
This is fairly unwieldy, but luckily going from $2$-cocycles to $2$-encoding tuples is easier.
\begin{proposition}[{\cite[p.~423]{cohom-abelian-crossed}}] \label{prop:cocycle-to-tuple}
    Let $A$ be a $G$-module, and let $c_2\in Z^2(G,A)$ be a $2$-cocycle. For each $i,j\in\{1,\ldots,m\}$, set
    \[\alpha_i\coloneqq\sum_{k=0}^{n_i}c_2\left(\sigma_i^k,\sigma_i\right)\qquad\text{and}\qquad\beta_{ij}\coloneqq c_2(\sigma_i,\sigma_j)-c_2(\sigma_j,\sigma_i).\]
    Then $(\alpha,\beta)\in A^m\times A^{m\times m}$ is a $2$-encoding tuple and represents the cohomology class of $c_2$ via the map in \Cref{prop:tuple-to-cocycle}.
\end{proposition}
\begin{remark}
    Roughly speaking, \Cref{prop:tuple-to-cocycle,prop:cocycle-to-tuple} are constructing quasi-isomorphisms between the standard free resolution of $\ZZ$ (which gives rise to cocycles) and another arising from writing $G$ as a product of cyclic groups. This cleanly explains where these maps come from and why they should be inverse.
\end{remark}

\subsubsection{Explicit Inflation--Restriction}

In this subsubsection, we make the Inflation--Restriction exact sequence explicit in dimension $2$. The results are similar to \cite[Section~2]{explicit-fund-classes}, but we will need to be more general. Roughly speaking, the proof is by making the degree-$1$ case explicit and then dimension-shifting with the short exact sequence
\begin{equation}
    0\to A\to\op{Hom}_\ZZ(\ZZ[G],A)\to\op{Hom}_\ZZ(I_G,A)\to0. \label{eq:the-dim-shift-ses}
\end{equation}
Here is the dimension-shifting.
\begin{lemma}[{\cite[Lemma~2.1]{explicit-fund-classes}}] \label{lem:explicitdimshift}
	Let $G$ be a group. Given a $G$-module $A$, the maps
	\begin{align*}
		\delta\colon Z^1(G,\op{Hom}_\ZZ(I_G,A))&\to Z^2(G,A) \\
		c_1&\mapsto\left[(g,g')\mapsto g\cdot c_1(g')(g^{-1}-1)\right] \\
		\left[g\mapsto\big((g'-1)\mapsto g'\cdot c_2((g')^{-1},g)\big)\right]&\mapsfrom c_2
	\end{align*}
	are inverse homomorphisms descending to the $\delta$-isomorphism $H^1(H,\op{Hom}_\ZZ(I_G,A))\cong H^2(H,A)$ of \eqref{eq:the-dim-shift-ses}.
\end{lemma}
\begin{proof}
    Computing the forward map $\delta$ is a matter of tracking through the $\delta$-morphism induced by \eqref{eq:the-dim-shift-ses}. It remains to show the last sentence. Given $u\in Z^2(H,A)$, define $c_u\in C^1(G,\op{Hom}_\ZZ(I_G,A))$ by
	\[c_u(g)(g'-1)=g'\cdot u\left((g')^{-1},g\right).\]
	One can verify directly that $c_u$ is a $1$-cocycle using the $2$-cocycle condition on $u$. Similarly, it is a direct computation that these maps are inverse.
\end{proof}
And here is our result. It is similar to \cite[Lemma~2.3]{explicit-fund-classes}.
\begin{proposition} \label{prop:explicitresinf}
	Let $G$ be a group with normal subgroup $H\subseteq G$. Fix a $G$-module $A$ with $H^1(H,A)=0$, and find some $c_2\in Z^2(G,A)$ such that $\op{Res}^G_H([c_2])=0$ in $H^2(H,A)$.
    \begin{enumerate}[label=(\alph*)]
        \item There exist elements $b_g\in A$ for each $g\in G$ such that $b_e=0$ and
        \[c_2(g,h)=b_g+gb_h-b_{gh}\]
        for each $g\in G$ and $h\in H$.
        \item There exist elements $\eta_g\in A$ for each $g\in G$ such that
        \[(h-1)\eta_g=c_2(h,g)+(b_{hg}-hb_g-b_h)\]
        for each $g\in G$ and $h\in H$. In fact, we may define the $\eta_g$ so that $\eta_g$ only depends on the coset $gH$.
        \item The $2$-cochain $u\in C^2(G,A)$ defined by
        \[u(g,g')\coloneqq c_2(g,g')+(b_{gg'}-gb_{g'}-b_g)+(\eta_{gg'}-g\eta_{g'}-\eta_g)\]
        actually defines a $2$-cocycle in $Z^2\left(G/H,A^H\right)$ which inflates to the class of $c$ in $H^2(G,A)$.
    \end{enumerate}
\end{proposition}
\begin{proof}
    These proofs are essentially direct computations. Define $\delta$ and $\delta^{-1}$ as in \Cref{lem:explicitdimshift}.
    \begin{enumerate}[label=(\alph*)]
        \item Note that $\op{Res}\delta^{-1}([c_2])=\delta^{-1}\op{Res}([c_2])=[0]$ in $H^1(H,\op{Hom}_\ZZ(I_G,A))$, so we can find $b_0\in\op{Hom}_\ZZ(I_G,A)$ such that $\op{Res}\delta^{-1}c_2=db_0$. To understand this condition, we note that we are asking for
        \[g^{-1}\cdot c_2(g,h)=\left(\delta^{-1}c_2\right)(h)\left(g^{-1}-1\right)=(db_0)(h)\left(g^{-1}-1\right)\]
        for any $g\in G$ and $h\in H$. Setting $b_g\coloneqq-gb_0\left(g^{-1}-1\right)$, we see that we are asking for elements $b_g\in A$ for each $g\in G$ such that $b_e=0$ and
        \[c_2(g,h)=b_g+gb_h-b_{gh}\]
        for each $g\in G$ and $h\in H$.
        
        \item Set $c_1\coloneqq\left(\delta^{-1}c_2-db_0\right)\in Z^1(G,\op{Hom}_\ZZ(I_G,A))$. One can expand out the definitions to show
        \[g\cdot c_1(g')\left(g^{-1}-1\right)=c_2(g,g')+(b_{gg'}-gb_{g'}-b_g).\]
        Now, by construction, $c_1$ vanishes on $H$, so the $1$-cocycle condition implies $c_1(gh)=c_1(g)$ for each $g\in G$ and $h\in H$, so the function $c_1(g)\in\op{Hom}_\ZZ(I_G,A)$ is only defined up to coset $gH$. Similarly, for each $g\in G$, we see $c_1$ vanishing on $H$ implies $hc_1(g)=c_1(g)$, so we can check that the $1$-cochain $c_g\in C^1(H,A)$ defined by
        \[c_g(h)\coloneqq c_1(g)(1-h)=h\cdot c_1(g)\left(h^{-1}-1\right)=c_2(h,g)+(b_{hg}-hb_g-b_h)\]
        is a $1$-cocycle in $Z^1(H,A)$. Because $H^1(H,A)=0$, we conclude that $c_g$ is a $1$-coboundary, and the result follows.

        \item It's enough to check that $u$ is a well-defined $2$-cochain in $C^2\left(G/H,A^H\right)$: indeed, $u$ is just $c_2$ with some added $2$-coboundaries, so it is a $2$-cocycle and will end up inflating properly to the class of $c$ in $H^2(G,A)$.
        
        Fix $g,g'\in G$ and $h\in H$. We begin by showing that $u$ is well-defined up to coset in $G/H$. Note
        \[u(g,g')=g\cdot c_1(g')\left(g^{-1}-1\right)+(\eta_{gg'}-g\eta_{g'}-\eta_g),\]
        so $u(g,g'h)=u(g,g')$ because $c_1$ and the $\eta$s are well-defined up to coset in $G/H$. Additionally, we can see $u(gh,g')=u(g,g')$ because $h\eta_{g'}=c_1(g')(1-h)+\eta_{g'}$.

        It remains to show $h\cdot u(g,g')=u(g,g')$. This follows by summing the identities
        \begin{align*}
            (h-1)g\cdot c_1(g')\left(g^{-1}-1\right) &= g\cdot c_1(g')\left(g^{-1}h-g^{-1}hg-g^{-1}-1\right) \\
            (h-1)(\eta_{gg'}-\eta_g) &= g\cdot c_1(g')\left(g^{-1}-g^{-1}h\right) \\
            (h-1)\cdot-g\eta_{g'} &= g\cdot c_1(g')\left(g^{-1}hg-1\right),
        \end{align*}
        each of which follow by direct expansion with the definitions.
        \qedhere
    \end{enumerate}
\end{proof}

\subsection{Computation} \label{subsec:local-computation}
In this subsection we record the details of the computation.

\subsubsection{Set-Up} \label{subsubsec:set-up}
Fix a finite abelian extension of local fields $L/K$. 
We now build the following tower of fields.
\begin{equation}
    \begin{tikzcd}
    	LM \\
    	L & M \\
    	{K_{\pi}} & {L\cap M} \\
    	& K
    	\arrow["{\text{ram}}", no head, from=4-2, to=3-1]
    	\arrow["{\text{unr}}", no head, from=3-1, to=2-1]
    	\arrow["{\text{unr}}"', no head, from=4-2, to=3-2]
    	\arrow["{\text{ram}}"', no head, from=3-2, to=2-1]
    	\arrow["{\text{unr}}"', no head, from=3-2, to=2-2]
    	\arrow["{\text{unr}}", no head, from=2-1, to=1-1]
    	\arrow["{\text{ram}}"', no head, from=2-2, to=1-1]
    \end{tikzcd} \label{eq:set-up-tower}
\end{equation}
Unramified extensions have been labeled by ``unr,'' and totally ramified extensions have been labeled by ``ram.'' Here is the construction.
\begin{itemize}
    \item Set $n\coloneqq[L:K]$, and let $M$ be the unramified extension of $K$ of degree $n$. Note $\op{Gal}(M/K)$ is cyclic of order $n$ generated by the Frobenius automorphism $\sigma_0$. Extend $\sigma_0$ to an automorphism in $\op{Gal}(LM/K)$.
    \item Set $K_\pi\coloneqq(LM)^{\langle\sigma_0\rangle}$ so that $LM/K_\pi$ is unramified. Note $\op{Gal}(K/K_\pi)$ is a a finite abelian group, so we can write
    \[\op{Gal}(K/K_\pi)=\bigoplus_{i=1}^m\langle\sigma_i\rangle\]
    for some $\sigma_i\in\op{Gal}(K_\pi/K)$ of order $n_i$.
    \item Note $K_\pi$ is totally ramified over $K$, and $M$ is unramified over $K$, so $K_\pi$ and $M$ are linearly disjoint. Thus,
    \[\op{Gal}(LM/K)=\op{Gal}(M/K)\oplus\op{Gal}(K_\pi/K)=\langle\sigma_0\rangle\oplus\bigoplus_{i=1}^m\langle\sigma_i\rangle.\]
    Extend each $\sigma_i$ up to an automorphism of $LM$ accordingly.
    \item Set $f\coloneqq[L\cap M:K]$. Note $L\cap M$ is the maximal unramified subextension of $L$, so $\op{Gal}((L\cap M)/K)$ is cyclic generated by $\sigma_0|_{L\cap M}$ and is of order $f$. Because of the isomorphism $\op{Gal}(LM/L)\cong\op{Gal}(M/(L\cap M))$, we see that $\op{Gal}(LM/L)$ is cyclic of order $n/f$ generated by $\sigma_0^f$. It follows $L=(LM)^{\langle\sigma_0^f\rangle}$, so $K_\pi\subseteq L$, and
    \[\op{Gal}(L/K)=\frac{\op{Gal}(LM/K)}{\op{Gal}(LM/L)}=\langle\sigma_0|_L\rangle\oplus\bigoplus_{i=1}^m\langle\sigma_i\rangle,\]
    where the restriction $\sigma_0|_L$ has order $n_0\coloneqq f$.
    \item Lastly, degree arguments or residue field computations can verify $LM/M$ and $L/(L\cap M)$ are all totally ramified extensions, though we will not need this.
\end{itemize}
For brevity, we will also set $G\coloneqq\op{Gal}(LM/K)$ and $H\coloneqq\op{Gal}(LM/L)$ so that $G/H=\op{Gal}(L/K)$. As an aside, we note that $\sigma_0|_L=\id_L$ when $f=1$, so we will not need to consider the generator $\sigma_0|_L$ in $\op{Gal}(L/K)$ when $L/K$ is totally ramified.

\subsubsection{Computing the Cocycle}
We are now ready to precisely state \Cref{thm:we-can-compute-fund-class} and provide its proof. Given a finite Galois extension $E/F$ of local fields, we will let $u_{E/F}$ denote the fundamental class in $H^2(\op{Gal}(E/F),E^\times)$ and $c_{E/F}$ a representative of $u_{E/F}$. Here is our statement.
\begin{theorem} \label{thm:fund-tuple}
	Fix everything as in \cref{subsubsec:set-up}. Letting $\pi$ be a uniformizer of $K$, there exists some $\gamma\in LM^\times$ such that $\op N_{LM/L}(\gamma)=\pi$ and elements $\eta_i\in LM^\times$ for $i\in\{0,1,\ldots,m\}$ such that
	\[\frac{\sigma_0^{f}\left(\eta_i\right)}{\eta_i}=\frac{\sigma_i(\gamma)}{\gamma}.\]
	Then the tuple given by
	\[\alpha_i\coloneqq\begin{cases}
        \gamma/\sigma_0^{(f)}(\eta_0) & \text{if }i=0, \\
		1/\sigma_i^{(n_i)}(\eta_i) & \text{if }i>0,
    \end{cases}\qquad\text{and}\qquad \beta_{ij}\coloneqq\frac{\sigma_j(\eta_i)}{\eta_i}\cdot\frac{\eta_j}{\sigma_i(\eta_j)},\]
	for each $i,j\in\{0,1,\ldots,m\}$ defines a $2$-encoding tuple representing the fundamental class $u_{L/K}$. If $L/K$ is totally ramified, then the same statement holds, ignoring the data given by index $0$.
\end{theorem}
\begin{proof}
    The idea is to chase the following diagram around.
    \begin{equation}
        \begin{tikzcd}
        	&& {H^2(\op{Gal}(M/K),M^\times)} \\
        	0 & {H^2(\op{Gal}(L/K),L^\times)} & {H^2(\op{Gal}(LM/K),LM^\times)} & {H^2(\op{Gal}(LM/L),LM^\times)}
        	\arrow["{\op{Res}}", from=2-3, to=2-4]
        	\arrow["{\op{Inf}}", from=2-2, to=2-3]
        	\arrow["{\op{Inf}}", from=1-3, to=2-3]
        	\arrow[from=2-1, to=2-2]
        \end{tikzcd} \label{eq:local-fund-diagram}
    \end{equation}
    Namely, using \cite[Lemma~2.7]{milne-cft}, we see
    \[\op{Inf}_{M/K}^{LM/K}u_{M/K}=[LM:M]u_{LM/K}=[LM:L]u_{LM/K}=\op{Inf}_{L/K}^{LM/K}u_{L/K}.\]
    Thus, we want to invert the rightward $\op{Inf}$ in \eqref{eq:local-fund-diagram}, for which we use \Cref{prop:explicitresinf} because the bottom row is an Inflation--Restriction exact sequence. Here are the steps. Throughout, we repeatedly write elements $g\in\op{Gal}(LM/K)$ as $\sigma_0^a\tau$ where $a\in\{0,1,\ldots,n-1\}$ and $\tau\in\op{Gal}(LM/M)$; note this representation is unique.
    \begin{enumerate}[label=(\alph*)]
        \item By \cite[p.~102]{milne-cft}, we may write
        \[c_{M/K}\left(\sigma_0^i,\sigma_0^j\right)\coloneqq \pi^{\floor{\frac{i+j}n}}=\begin{cases}
        	1 & \text{if }i+j<n, \\
        	\pi & \text{if }i+j\ge n,
        \end{cases}\]
        We claim $\pi\in\op N_{LM/L}(LM^\times)$. Indeed, let $\varpi$ be a uniformizer of $K_\pi$, and we note
        \[\op N_{LM/L}(\varpi)=\varpi^{[LM:L]}=\varpi^{[LM:M]}=\varpi^{[K_\pi:K]}\]
        has the same valuation as $\pi$. Thus, $\pi/\op N_{LM/L}(\varpi)\in\OO_L^\times$ and is thus in $\op N_{LM/L}(LM^\times)$ because $LM/L$ is unramified. It follows we can find $\gamma\in LM^\times$ such that $\op N_{LM/L}(\gamma)=\pi$.
    
        To define $b_g$, we write $g=\sigma_0^a\tau$ where $a\in\{0,1,\ldots,n-1\}$ and $\tau\in\op{Gal}(K_\pi/K)$. Then we define
        \[b_{\sigma_0^a\tau}\coloneqq\sigma_0^a\tau\Bigg(\prod_{i=1}^{\floor{a/f}}\sigma_0^{-fi}(\gamma)\Bigg).\]
        One can directly check that
        \[\frac{b_{\sigma_0^a\tau}\cdot\sigma_0^a\tau(b_{\sigma_0^{fb}})}{b_{\sigma_0^{a+fb}\tau}}=\sigma_0^a\tau\pi^{\floor{\frac{a+fb}n}}=(\op{Inf}c_{M/K})\left(\sigma_0^a\tau,\sigma_0^{fb}\right)\]
        for any $\sigma_0^a\tau\in G$ and $\sigma_0^{fb}\in H$. This is what we wanted.
        
        Now, we take a moment to write
        \begin{equation}
            c'(g,g')\coloneqq(\op{Inf}c_{M/K})(g,g')\cdot\frac{b_{gg'}}{gb_{g'}\cdot b_g}=\pi^{\floor{(a+a')/n}}\cdot\frac{\displaystyle gg'\Bigg(\prod_{i=1}^{\floor{[a+a']_n/f}}\sigma_0^{-fi}(\gamma)\bigg/\prod_{i=1}^{\floor{a'/f}}\sigma_0^{-fi}(\gamma)\Bigg)}{\displaystyle g\prod_{i=1}^{\floor{a/f}}\sigma_0^{-fi}(\gamma)} \label{eq:almost-local-fund}
        \end{equation}
        for $g=\sigma_0^a\tau$ and $g'=\sigma_0^{a'}\tau'$ in $G$. Here, $[a+a']_n$ refers to the remainder when $a+a'$ is divided by $n$.
    
        \item We discuss some of the $\eta$ terms. For brevity, define $\eta_i\coloneqq\eta_{\sigma_i}$ for each $\sigma_i$. Now, for $\tau\in\op{Gal}(K_\pi/K)$, we use \eqref{eq:almost-local-fund} to see $\eta_\tau$ satisfies
        \[\frac{\sigma_0^{fa}(\eta_\tau)}{\eta_\tau}=\sigma_0^{fa}\prod_{i=1}^{a}\frac{\tau\sigma_0^{-fi}(\gamma)}{\sigma_0^{-fi}(\gamma)}=\prod_{i=0}^{a-1}\frac{\tau\sigma_0^{fi}(\gamma)}{\sigma_0^{fi}(\gamma)}\]
        for each $\sigma_0^{fa}\in H$. By telescoping, it is enough for $\eta_\tau$ to satisfy the above equation at $a=1$, so we are asking for $\sigma_0^f(\eta_\tau)/\eta_\tau=\tau(\gamma)/\gamma$. Letting $\tau=\sigma_i$ for $i\in\{1,2,\ldots,m\}$ defines the necessary $\eta_i$ elements, and the argument at index $0$ is similar.
    
    
        \item We set
        \[c_{L/K}(g,g')\coloneqq c'(g,g')\cdot\frac{\eta_{gg'}}{g(\eta_{g'})\cdot \eta_g}\]
        to represent the fundamental class of $L/K$. We now compute the tuple for $c_{L/K}$ using \Cref{prop:cocycle-to-tuple}. For each $i,j\in\{1,2,\ldots,m\}$, we can use \eqref{eq:almost-local-fund} to compute
        \[\beta_{ij}=\frac{c'(\sigma_i,\sigma_j)}{c'(\sigma_j,\sigma_i)}\cdot\frac{\eta_{\sigma_i\sigma_j}/(\sigma_i(\eta_{\sigma_j})\cdot\eta_{\sigma_i})}{\eta_{\sigma_i\sigma_j}/(\sigma_j(\eta_{\sigma_i})\cdot\eta_{\sigma_j})}=\frac{\sigma_j(\eta_i)}{\eta_i}\cdot\frac{\eta_j}{\sigma_i(\eta_j)},\]
        and
        \[\alpha_i=\prod_{k=0}^{n_i-1}\left(c'\left(\sigma_i^k,\sigma_i\right)\cdot\frac{\eta_{\sigma_i^{k+1}}}{\sigma_i^k(\eta_{\sigma_i})\cdot\eta_{\sigma_i^k}}\right)=\sigma_i^{(n_i)}(\eta_i)^{-1}.\]
        When $L/K$ is not totally ramified so that $f>1$, the arguments at index $0$ are similar.
        \qedhere
    \end{enumerate}
\end{proof}
\begin{remark}
    \Cref{thm:fund-tuple} provides an improvement on Serre's approach (see \cite[Exercise~XIII.5.2]{local-fields} or \cite[Algorithm~2.18]{debeerst-cft-algorithms}) by working with a fixed finite unramified extension $M$ of $K$ instead of the full unramified closure $K^{\mathrm{unr}}$. Doing so makes an explicit algorithm easier computationally.
\end{remark}
\begin{remark} \label{rem:general-local-fund}
    Only steps (b) and (c) of the above proof use the fact that $\op{Gal}(L/K)$ is abelian. As such, we note \eqref{eq:almost-local-fund} provides a computation in the local fundamental class in full generality.
\end{remark}
\begin{remark}
    It is conceivable that a similar approach could work to compute global fundamental classes. For example, given a finite extension of number fields $L/K$, one can find a prime $p$ such that $\QQ(\zeta_p)\cap L=\QQ$ and $[L:K]\mid p-1$. Then $K(\zeta_p)$ can play the role of $M$ in the above discussion: with care, one can inflate a global fundamental class of $K(\zeta_p)/K$ up to $L(\zeta_p)/K$ and then restrict it back to $L/K$. The computations will be less explicit because one must work with id\'ele classes.
\end{remark}
\begin{remark}
    The tuples constructed in \Cref{thm:fund-tuple} behave in ``deeply ramified'' extensions. Namely, given a sequence of totally ramified extensions
    \[K=K_0\subseteq K_1\subseteq K_2\subseteq\cdots,\]
    let $M_i$ denote the unramified extension of $K$ of degree $[K_i:K]$. Then one can construct $\gamma_i\in(K_iM_i)^\times$ for each $i\in\NN$ with $\gamma_0=\pi$ and $\op N_{K_iM_i/K_iM_j}(\gamma_i)=\gamma_j$ for each $i\ge j$. These $\gamma_i$ allow construction of tuples which cohere with each other somewhat.
\end{remark}
As an application, we can compute the Artin reciprocity map. Given a finite Galois extension $E/F$ of local fields, we let $\theta_{E/F}\colon F^\times\to\op{Gal}(E/F)^{\mathrm{ab}}$ denote the Artin reciprocity map. Here is our result; it is essentially \cite[Theorem~1]{explicit-local-recip-dwork}.
\begin{corollary} \label{cor:explicit-local-recip}
    Fix notation as in \Cref{thm:fund-tuple}. Then $\theta_{L/K}(\op N_{L^{\langle\sigma_i\rangle}/K}(\alpha_i))=\sigma_i|_L$ for each $i$.
\end{corollary}
\begin{proof}
    Let $c_{L/K}$ be the $2$-cocycle coming from the $2$-encoding tuple of \Cref{thm:fund-tuple}, and for brevity, we set $L_i\coloneqq L^{\langle\sigma_i\rangle}$ so that $\op{Gal}(L_i/K)=\langle\sigma_i|_L\rangle$ has order, where by $\sigma_0$ we mean the restriction of $\sigma_0$ to $L$. The main point is that, by definition of $\alpha_i$, we have
    \[\alpha_i=\prod_{k=0}^{n_i-1}c_{L/K}\left(\sigma_i^k,\sigma_i\right)=(\op{Res}_{\langle\sigma_i|_L\rangle}u_{L/K})\cup[\sigma_i],\]
    where we used the cup-product computation of \cite[Proposition~I.5.8]{bonn-lectures}. However, \cite[Lemma~2.7]{milne-cft} promises $\op{Res}_{\langle\sigma_i|_L\rangle}u_{L/K}=u_{L/L_i}$, so $\alpha_i=\theta^{-1}_{L/L_i}(\sigma_i)$. To finish, we note $\theta_{L/K}\circ{\op N_{L_i/K}}=\theta_{L/L_i}$ by \cite[Lemma~3.2]{milne-cft}.
\end{proof}
\begin{remark}
    Because the $\sigma_i|_L$ generate $\op{Gal}(L/K)$, and $\theta_{L/K}$ is an isomorphism $K^\times/\op N_{L/K}(L^\times)\to\op{Gal}(L/K)$, we see that the elements
    \[\op N_{L^{\langle\sigma_i\rangle}/K}(\alpha_i)\]
    generate $K^\times/\op N_{L/K}(L^\times)$. Thus, \Cref{cor:explicit-local-recip} does determine the reciprocity map.
\end{remark}
\begin{remark} \label{rem:cannot-get-betas-from-artin}
    \Cref{cor:explicit-local-recip} shows us that the $\beta$s of the $2$-encoding tuple cannot be recovered from merely understanding the Artin reciprocity map: the Artin reciprocity map is entirely determined by the $\alpha$s!
\end{remark}

\subsubsection{Tame Ramification}

In this section, we work through \autoref{thm:fund-tuple} for some tamely ramified extensions in order to prove \Cref{prop:tame-tuple}.
\begin{proof}[Proof of \Cref{prop:tame-tuple}]
    To show that we are within the set-up described in \cref{subsubsec:set-up}, we will show that the extension $L/K$ fills in \eqref{eq:set-up-tower} as the following tower of fields.
    \begin{equation}
        \begin{tikzcd}
        	{K(\gamma,\zeta_{q^{ef}-1})} \\
        	{L} & {K(\zeta_{q^{ef}-1})} \\
        	{K(\gamma)} & {K(\zeta_{q^f-1})} \\
        	& {K}
        	\arrow[no head, from=4-2, to=3-2]
        	\arrow[no head, from=3-2, to=2-2]
        	\arrow[no head, from=4-2, to=3-1]
        	\arrow[no head, from=3-1, to=2-1]
        	\arrow[no head, from=3-2, to=2-1]
        	\arrow[no head, from=2-1, to=1-1]
        	\arrow[no head, from=2-2, to=1-1]
        \end{tikzcd} \label{eq:tame-ram-tower}
    \end{equation}
    Indeed, we now construct the required elements.
    \begin{itemize}
        \item Note that the extension $L^{\langle\sigma_0\rangle}/K$ is totally ramified with ramification index $e$, so it is tamely ramified also. Thus, \cite[Proposition~II.5.12]{lang-alg-nt} promises some element $\gamma\in L^{\langle\sigma_0\rangle}/K$ such that $\gamma^e=\pi$. Notably, $K(\gamma)$ has degree $e$ over $K$, so we conclude $K(\gamma)=L^{\langle\sigma_0\rangle}$. This explains the $K(\gamma)$ in \eqref{eq:tame-ram-tower}.
        \item Note $L$ has residue field $\FF_{q^f}$ by definition of $f$. The polynomial $x^{q^f}-1=0$ splits into linear factors over $\FF_{q^f}$, so we can use Hensel's lemma to lift any root up to $L$. Thus, we can find $\zeta_{q^f-1}\in L$; this also explains how we find $\zeta_{q^{ef}-1}$ in the unramified extension of degree $[L:K]=ef$ over $K$ in \eqref{eq:tame-ram-tower}.

        Conversely, we note that $K(\zeta_{q^f-1})\subseteq L$ has residue field containing $\FF_{q^f}$, so we must have $[K(\zeta_{q^f-1}):K]=f$, meaning that $K(\zeta_{q^f-1})$ is the maximal unramified subextension of $L/K$. A similar argument shows that $K(\zeta_{q^{ef}-1})$ is the unramified extension of degree $[L:K]=ef$ over $K$.
    \end{itemize}
    We now take a moment to see $\sigma_0$ can be extended to $K(\gamma,\zeta_{q^{ef}-1})/K$ by fixing $\gamma$ but sending $\zeta_{q^{ef}-1}\mapsto\zeta_{q^{ef}-1}^q$. Similarly, $\sigma_1$ is extended by fixing $\zeta_{q^{ef}-1}$ but sending $\gamma\mapsto\zeta_e\gamma$. To tie up loose ends, we note that the Galois conjugates of $\gamma$ in $K(\gamma)$ take the form $\zeta_e^i\gamma=\sigma_1^i(\gamma)$, so $\sigma_1$ does indeed generate $\op{Gal}(K(\gamma)/K)$.
    
    We now finish the proof. Using $\pi$ as our uniformizer, we may use our found $\gamma$ on \Cref{thm:fund-tuple}. For example, $\gamma\in K(\gamma)$ is fixed by $\sigma_0$, so we may set $\eta_0=1$, yielding $\alpha_0=\gamma$. It remains to compute $\eta_1$. Observe that
    \[b\coloneqq\frac{q^{ef}-1}{e\left(q^f-1\right)}\]
    is an integer because $b=\frac1e\sum_{k=0}^{e-1}q^{kf}$ and $e\mid q-1$. As such, we define $\eta_1\coloneqq\zeta_{q^{ef}-1}^b$. Indeed,
    \[\frac{\sigma_0^f(\eta_1)}{\eta_1}=\frac{\zeta_{bq^{ef}-1}^{q^f}}{\zeta_{q^{ef}-1}^b}=\zeta_{q^{ef}-1}^{b\left(q^f-1\right)}=\zeta_e.\]
    To finish, one can use the formulae provided in \Cref{thm:fund-tuple} to directly compute $\alpha_1$ and the $\beta_{ij}$.
\end{proof}
\begin{remark}
    Roughly speaking, we require $e\mid q-1$ to ensure a primitive $e$th root of unity lives in $K$ so that the extension $K(\gamma)/K$ is cyclic. It might appear that the proof can go through with the assumption $\zeta_e\in K(\gamma)$ instead of $e\mid q-1$, but these conditions are in fact equivalent: we showed the reverse direction in the proof, and if $\zeta_e\in K(\gamma)$, then the residue field extension $\FF_q$ of $K(\gamma)$ has an element of order $e$, so it follows $e\mid q-1$.
\end{remark}
\begin{example} \label{ex:basic-ramified-artin-map}
    We use \Cref{prop:tame-tuple} to verify Artin reciprocity for $\QQ_p(\zeta_p)/\QQ_p$; note that this extension satisfies the hypotheses. Let $x\in\FF_p^\times$ be a generator so that the automorphism $\sigma_1\colon\zeta_p\mapsto\zeta_p^x$ generates $\op{Gal}(\QQ_p(\zeta_p)/\QQ_p)$. Letting $\gamma$ be an element such that $\gamma^{p-1}=p$, we note that $\zeta_{p-1}\coloneqq\sigma_1(\gamma)/\gamma$ is a $(p-1)$st root of unity; we claim $\zeta_{p-1}\equiv x\pmod p$. By the uniqueness of the Hensel-lifting constructing our $(p-1)$st roots of unity, it is enough to show
    \[\left|\frac{\sigma_1(\gamma)}{\gamma}-x\right|<1,\]
    where $|\cdot|$ is the usual absolute value. Let $\varpi\coloneqq(\zeta_p-1)$ be a uniformizer of $\QQ_p(\zeta_p)$. Note $|\gamma|=p^{-1/(p-1)}=|\varpi|$, so $\gamma$ is also a uniformizer of $\QQ_p(\zeta_p)$. Thus, $\gamma/\varpi$ reduces to a unit in $\OO_{\QQ_p(\zeta_p)}/\varpi\OO_{\QQ_p(\zeta_p)}\cong\FF_p$, so there is $c\in\ZZ\setminus p\ZZ$ such that $\gamma/\varpi\equiv c\pmod{\varpi}$. It follows
    \[\frac{\sigma_1(\gamma)}{\gamma}\equiv\frac{c\cdot\sigma_1(\varpi)}{c\cdot\varpi}=\frac{\zeta_p^x-1}{\zeta_p-1}=1+\zeta_p+\cdots+\zeta_p^{x-1}\equiv x\pmod{\varpi},\]
    so we have shown the claim. Finishing up, we see
    \[\alpha_1=\zeta_{p-1}^{-1}\equiv x^{-1}\pmod p.\]
    Combining with \Cref{cor:explicit-local-recip}, we see that the inverse Artin map $\theta^{-1}_{\QQ_p(\zeta_p)/\QQ_p}$ sends $\sigma_1\colon\zeta_p\mapsto\zeta_p^x$ to $x^{-1}\pmod p$, as predicted by Lubin--Tate theory \cite[Theorem~3.2]{local-cft}.
\end{example}

\section{The Tate Canonical Class} \label{sec:global}

In this section, we execute \Cref{thm:we-can-tate-class} and compute the Tate canonical class for the extension $\QQ(\zeta_{p^\nu})/\QQ$, where $p^\nu$ is an odd prime-power.

\subsection{Construction} \label{subsec:tate-class-defi}
The goal of this subsection is to define the Tate canonical class so that we may compute in the subsequent subsections. We will directly follow \cite{tate-torus} and \cite{kottwitz}. Given a global field $F$, let $V_F$ denote the set of places of $F$. 


Now, fix a finite Galois extension of global fields $L/K$ with Galois group $G\coloneqq\op{Gal}(L/K)$. For later use, we will also let $G_v\subseteq G$ denote the decomposition group of a place $v\in V_L$.
Now, we build two short exact sequences, as described in \autoref{sec:intro}. To begin, we note that the augmentation map $\ZZ[V_L]\onto\ZZ$ induces the short exact sequence
\[0\to\ZZ[V_L]_0\to\ZZ[V_L]\to\ZZ\to0\label{eq:sesx}\tag{$X$}\]
where $\ZZ[V_L]_0$ is the kernel of $\ZZ[V_L]\onto\ZZ$. We also have the short exact sequence
\[0\to L^\times\to\AA_{L}^\times\to\AA_L^\times/L^\times\to0\tag{$A$}\label{eq:sesa}\]
where $L^\times\into\AA_L^\times$ is the diagonal embedding.

We now construct three different cohomology classes, $\alpha_1$, $\alpha_2$, and $\alpha_3$. Roughly speaking, $\alpha_1$ will be constructed from global class field theory, $\alpha_2$ from local class field theory, and $\alpha_3$ will measure their compatibility. To begin, we let $\alpha_1\in\widehat H^2\left(G,\op{Hom}_\ZZ(\ZZ,\AA_L^\times/L^\times)\right)$ denote the global fundamental class.

We now construct $\alpha_2\in H^2(G,\op{Hom}_\ZZ(\ZZ[V_L],\AA^\times_{L}))$ from local class field theory. We need the following lemma; here, for each place $u\in V_K$, we fix a place $v(u)\in V_L$ above $u$.
\begin{lemma}[{\cite[p.~714]{tate-torus}}] \label{lem:magicaltate}
    Fix everything as above, and let $M$ be a $G$-module. Then, for any $i\in\ZZ$, the morphisms
    \[\widehat H^i(G,\op{Hom}_\ZZ(\ZZ[V_L],M))\xrightarrow{\op{Res}}\widehat H^i(G_{v},\op{Hom}_\ZZ(\ZZ[V_L],M))\xrightarrow{\op{eval}_v}\widehat H^i(G_{v},M)\]
    glue into an isomorphism
    \[\widehat H^i(G,\op{Hom}_\ZZ(\ZZ[V_L],M))\cong\prod_{u\in S_K}\widehat H^i(G_{v(u)},M),\]
    where the product is over places $u\in V_K$ with a chosen place $v(u)\in V_L$ above $u$.
\end{lemma}
Thus, to specify $\alpha_2\in H^2(G,\op{Hom}_\ZZ(\ZZ[V_L],\AA_L^\times))$, it is enough to specify a set of classes
\[\alpha_2(u)\in H^2\left(G_{v(u)},\AA_{L}^\times\right)\]
for each $u\in S_K$. To do so, we note that $G_{v(u)}=\op{Gal}(L_{v(u)}/K_u)$, so we use the natural embedding $i_{v(u)}\colon L_{v(u)}\into\AA_{L}^\times$ to set
\[\alpha_2(u)\coloneqq i_{v(u)}\big(\alpha(L_{v(u)}/K_u)\big),\]
where $\alpha(L_{v(u)}/K_u)\in H^2\big(G_{v(u)},L_{v(u)}^\times\big)$ is the local fundamental class.

Lastly, we construct $\alpha_3\in H^2(G,\op{Hom}_\ZZ(\ZZ[V_L]_0,\OO_{L,V_L}^\times))$. Roughly speaking, we note that the morphism of short exact sequences from \autoref{eq:sesx} and \autoref{eq:sesa} can be specified by commuting morphisms $\ZZ[V_L]\to\AA_L^\times$ and $\ZZ\to\AA_L^\times/L^\times$, inducing the last arrow as follows.
\[\begin{tikzcd}
    0 & {\ZZ[V_L]_0} & {\ZZ[V_L]} & \ZZ & 0 \\
    0 & {L^\times} & {\AA_{L}^\times} & {\AA_L^\times/L^\times} & 0
    \arrow[from=1-1, to=1-2]
    \arrow[from=1-2, to=1-3]
    \arrow[from=1-3, to=1-4]
    \arrow[from=1-4, to=1-5]
    \arrow[from=2-1, to=2-2]
    \arrow[from=2-2, to=2-3]
    \arrow[from=2-3, to=2-4]
    \arrow[from=2-4, to=2-5]
    \arrow[dashed, from=1-2, to=2-2]
    \arrow[from=1-3, to=2-3]
    \arrow[from=1-4, to=2-4]
\end{tikzcd}\]
Intuitively, this should let us specify a cohomology class $\alpha_3\in H^2\big(G,\op{Hom}_\ZZ(\ZZ[V_L]_0,L^\times)\big)$ from melding together $\alpha_1$ and $\alpha_2$.

To rigorize this, we let $\op{Hom}_\ZZ(X,A)$ denote the group of morphisms of short exact sequences from \autoref{eq:sesx} to \autoref{eq:sesa}, and we let $\pi_3$, $\pi_2$, and $\pi_1$ denote the projections from $\op{Hom}_\ZZ(X,A)$ to $\op{Hom}_\ZZ(\ZZ[V_L]_0,L^\times)$, to $\op{Hom}_\ZZ(\ZZ[V_L],\AA_{S_L}^\times)$, and to $\op{Hom}_\ZZ(\ZZ,\AA_L^\times/L^\times)$ respectively. Extending the above argument is able to prove the following.
\begin{lemma}[{\cite[p.~716]{tate-torus}, \cite[Lemma~6.3]{kottwitz}}] \label{lem:constructtatealpha}
    Fix everything as above. Then
    \[\begin{tikzcd}
        {H^2(G,\op{Hom}_\ZZ(X,A))} & {H^2\big(G,\op{Hom}_\ZZ(\ZZ[V_L],\AA_{L}^\times)\big)} \\
        {H^2(G,\op{Hom}_\ZZ(\ZZ,\AA_L^\times/L^\times))} & {H^2\big(G,\op{Hom}_\ZZ(\ZZ[V_L],\AA_L^\times/L^\times)\big)}
        \arrow["{\pi_2}", from=1-1, to=1-2]
        \arrow["{\pi_1}"', from=1-1, to=2-1]
        \arrow[from=1-2, to=2-2]
        \arrow[from=2-1, to=2-2]
    \end{tikzcd}\]
    is a pull-back square.
\end{lemma} 
To finish, one can check that $\alpha_1$ and $\alpha_2$ have the same image in $H^2\big(G,\op{Hom}_\ZZ(\ZZ[V_L],\AA_{L}^\times/L^\times)\big)$ so that \Cref{lem:constructtatealpha} promises us $\alpha\in H^2(G,\op{Hom}_\ZZ(X,A))$ such that $\pi_2(\alpha)=\alpha_2$ and $\pi_1(\alpha)=\alpha_1$. These together let us construct $\alpha_3\coloneqq\pi_3(\alpha)\in H^2\big(G,\op{Hom}_\ZZ(\ZZ[V_L]_0,\AA_L^\times/L^\times)\big)$, which is the Tate canonical class.

\subsection{Idea of Computation}
We continue in the context of \cref{subsec:tate-class-defi}. We will compute the Tate canonical class $\alpha_3$ of the global extension $\QQ(\zeta_{p^\nu})/\QQ$, where $p$ is an odd prime, and $\nu$ is a positive integer. The idea is to follow the construction provided in \cref{subsec:tate-class-defi}. To begin, note that $(p)$ totally ramifies in $\QQ(\zeta_{p^\nu})$, so its decomposition group is the full Galois group. Thus,
we can compute a representative $c_1$ for $\alpha_1\in H^2(G,\AA_L^\times/L^\times)$ by using a local fundamental class.

For $\alpha_2$, it is relatively straightforward to give some representative $c_2$ for $\alpha_2\in H^2\big(G,\op{Hom}_\ZZ(\ZZ[V_L],\AA_{L}^\times)\big)$ by using, say, \Cref{thm:fund-tuple} to get representatives for the local fundamental classes and then inverting the isomorphism of \Cref{lem:magicaltate}.
In fact, in our toy case this is easier because all relevant extensions are cyclic.

However, the main difficulty here is in melding $c_2$ and $c_1$ to make the diagram
\[\begin{tikzcd}
    {\ZZ[V_L]} & \ZZ \\
    {\AA_{L}^\times} & {\AA_L^\times/L^\times}
    \arrow["", from=1-1, to=1-2]
    \arrow["{c_1(g,g')}", from=1-2, to=2-2]
    \arrow["{c_2(g,g')}"', from=1-1, to=2-1]
    \arrow["", from=2-1, to=2-2]
\end{tikzcd}\]
commute for all $g,g'\in G$. Thus, in \Cref{subsec:compute-a2}, we will fix $c_1$ and then choose a special $c_2$ to make this diagram commute. From here, we induce $c_3(g,g')$ for $g,g'\in\op{Gal}(L/K)$ by filling in the leftmost arrow in the following morphism of short exact sequences.
\[\begin{tikzcd}
    0 & {\ZZ[V_L]_0} & {\ZZ[V_L]} & \ZZ & 0 \\
    0 & {L^\times} & {\AA_{L}^\times} & {\AA_L^\times/L^\times} & 0
    \arrow[from=1-3, to=1-4]
    \arrow["{c_1(g,g')}", from=1-4, to=2-4]
    \arrow["{c_2(g,g')}", from=1-3, to=2-3]
    \arrow[from=1-2, to=1-3]
    \arrow[from=1-1, to=1-2]
    \arrow[from=2-1, to=2-2]
    \arrow[from=2-2, to=2-3]
    \arrow[from=1-4, to=1-5]
    \arrow[from=2-4, to=2-5]
    \arrow[from=2-3, to=2-4]
    \arrow["{c_1(g,g')}", dashed, from=1-2, to=2-2]
\end{tikzcd}\]
We do this in \cref{subsec:compute-a3} to compute a representative $c_3$ for $\alpha_3$.

\subsection{Setting Variables} \label{subsec:compute-a2}
We continue in the context of \autoref{subsec:tate-class-defi}, in the case of $K\coloneqq\QQ$ and $L\coloneqq\QQ(\zeta_{p^\nu})$. For brevity, set $\zeta\coloneqq\zeta_{p^\nu}$, and let $\mf P\coloneqq(\zeta-1)$ be the prime of $L$ above $(p)$. The goal of the computation is to
to be able to write down a $2$-cocycle in $Z^2\big(G,\op{Hom}_\ZZ(\ZZ[V_L],\AA_{L}^\times)\big)$ representing $\alpha_2$ compatible with a choice of global fundamental class.

\subsubsection{Extracting Elements}
We are going to choose our local fundamental class representatives to be compatible with a choice of global fundamental class for $L/K$. However, this will require extracting certain magical elements of $L^\times$, so we will go ahead and extract these before getting into the computation. Here are the elements we will need.
\begin{enumerate}
    \item We choose a generator $x$ of $\left(\ZZ/p^\nu\ZZ\right)^\times$ so that $\sigma\colon\zeta\mapsto\zeta^x$ is a generator of $G$. \label{item:choicegen}
    \item For any infinite place $v\mid\infty$ and subgroup $H\subseteq G$ containing $G_v$, we construct $\xi_{v,H}\in L^H$ such that
    \[\xi_{v,H}\equiv i_\mf P(x)^{|H|/2}\cdot\prod_{w\in Hv}i_w(-1)\pmod{N_H\AA_{L}^\times}.\] \label{item:choicexi}
    \item For each subgroup $H\subseteq G$ and ideal class $c\in\op{Cl}L^H$, we find a prime $L^H$-ideal $\mf r_{H,c}$ representing $c$ while splitting completely in $L$. \label{item:choicecheb}
    \item Given a finite unramified place $u\in V_K$ under $v(u)\in V_L$, set $\mf q$ to be the prime $L^{G_{v(u)}}$-ideal below $v(u)$. Finding the correct $\mf r_{H,[\mf q^{-1}]}$, we need to select a generator $\varpi_u$ for $\mf q\mf r_{G_{v(u)},[\mf q^{-1}]}$. \label{item:choiceprincgen}
\end{enumerate}
All of these elements above will be found somewhat non-constructively.

We now construct \eqref{item:choicegen}--\eqref{item:choiceprincgen}, in order. For \eqref{item:choicegen}, there is nothing to say, so we move on to \eqref{item:choicexi}.
By \cite[p.~102]{milne-cft},
\[c_p(\sigma^i,\sigma^j)=x^{-\floor{\frac{i+j}n}}\]
is a $2$-cocycle in $Z^2(G,L_\mf P^\times)$ representing the local fundamental class of $L_\mf P/K_{(p)}$. Passing $c_\mf P$ through $L_\mf P^\times\into\AA_L^\times\onto\AA_L^\times/L^\times$, we see that
\[i_\mf Pc_p(\sigma^i,\sigma^j)=i_\mf P(x)^{-\floor{\frac{i+j}n}}\]
has cohomology class of global invariant $1/n$ and therefore represents the global fundamental class $u_{L/K}\in H^2(G,\AA_L^\times/L^\times)$. We now begin constructing our $\xi$ elements. We begin at the fixed infinite place $v(\infty)$. Set $\tau\coloneqq\sigma^{n/2}$ to be the conjugation automorphism on $L$.
\begin{lemma}
    Let $v\coloneqq v(\infty)$ be our chosen infinite place so that $G_v=\{1,\tau\}$. Then there exists $\xi_\infty\in L^{\langle\tau\rangle}$ such that
    \[\xi_\infty\equiv i_v(-1)\cdot i_\mf P(x)\pmod{N_{\langle\tau\rangle}\AA_{L}^\times}.\]
\end{lemma}
\begin{proof}
    We can represent the local fundamental class of $L_v/K_\infty$ by the $2$-cocycle
    \[c_v(\tau^i,\tau^j)=(-1)^{\floor{\frac{i+j}2}}.\]
    Mapping $c_v$ to $\AA_L^\times/L^\times$, we see $[i_v(c_v)]\in H^2(G_v,\AA_L^\times/L^\times)$
    has global invariant $1/2$ and therefore should live in the same cohomology class as $\op{Res}_{G_v}i_\mf P(c_\mf P)$. In particular, using $[\tau]\in\widehat H^{-2}(G_v,\ZZ)$
    \[[i_v(-1)]=[i_v(c_v)]\cup[\tau]=[\op{Res}_{G_v}i_\mf P(c_\mf P)]\cup[\tau]=[i_\mf P(x)]^{-1}\]
    as elements in $\widehat H^0(G_v,\AA_L^\times/L^\times)$.
    Now, this group is $\AA_{L}^\times/L^\times$ modded out by $N_{G_v}\AA_{L}^\times$, so we can unwind the equivalence relations as promising some $\xi_\infty\in L^\times$ such that
    \[\xi_\infty\equiv i_v(-1)\cdot i_\mf P(x)\pmod{N_{G_v}\AA_{L}^\times}.\]
    Lastly, $\xi_\infty$ is fixed by $\tau$ because $i_v(-1)$ and $i_\mf P(x)$ are both fixed by $\tau$.
\end{proof}
\begin{remark}
    For certain primes $p$, one can choose $\xi_\infty$ from the circulant units of $\QQ(\zeta_p)$, making $\xi_\infty$ effectively computable with some linear algebra. In particular, \cite[Proposition~1]{signature-ranks-dummit} tells us that this is possible whenever $|\op{Cl}H|$ is odd; this fails first for $\QQ(\zeta_{29})$. In \cite{signature-ranks-david}, it is conjectured that this is possible for infinitely many primes.
\end{remark}
Extending to our general infinite places, we note that any infinite place takes the form $gv(\infty)$ for some $g\in G$, where we note
\[g\xi_\infty\equiv i_{gv}(-1)\cdot i_\mf P(x)\pmod{N_{G_{gv(\infty)}}\AA_L^\times}\]
because $G_{v(\infty)}=G_{gv(\infty)}$. Thus, we set $\xi_{gv(\infty)}\coloneqq g\xi_\infty\in L^{\langle\tau\rangle}$. Because $\xi_\infty$ is fixed by $\tau$, we see $\xi_{gv(\infty)}$ does not depend on the precise choice of $g\in G$.

To finish our construction of the $\xi$ elements, we need to account for subgroups $H\subseteq G$, which is the content of the following lemma.
\begin{lemma}
    Fix everything as above. Picking any infinite place $v\mid\infty$ and subgroup $H\subseteq G$ containing $\tau$, the element
    \[\xi_{v,H}\coloneqq\prod_{g\langle\tau\rangle\in H/\langle\tau\rangle}g\xi_v\]
    has
    \[\xi_{v,H}\in L^H\qquad\text{and}\qquad\xi_{v,H}\equiv i_\mf P(x)^{|H|/2}\cdot\prod_{w\in Hv}i_w(-1)\pmod{N_H\AA_{L}^\times}.\]
\end{lemma}
Again, we note that the choice of coset representatives in $H/\langle\tau\rangle$ does not matter to define $\xi_{v,H}$ because each $\xi_v$ is fixed by $\tau$.
\begin{proof}
    By construction,
    \[\xi_v=i_\mf P(x)\cdot i_v(-1)\cdot N_{\langle\tau\rangle}(a)\]
    for some $a\in\AA_{L,S_L}^\times$. Now, direct expansion shows
    \[\xi_{v,H}=i_\mf P(x)^{|H|/2}\cdot\Bigg(\prod_{w\in Hv}^mi_w(-1)\Bigg)\cdot N_H(a).\]
    This shows the last claim. Lastly,
    we see $\xi_{v,H}$ is fixed by $H$ because each factor above is fixed by $H$.
\end{proof}
Next we turn to our finite unramified places, the subject of \eqref{item:choicecheb} and \eqref{item:choiceprincgen}. There is nothing to say for the construction of \eqref{item:choiceprincgen}, so here is \eqref{item:choicecheb}.
\begin{lemma} \label{lem:magicalprimes}
    Fix everything as above. For each subgroup $H\subseteq G$ and ideal class $c\in\op{Cl}L^H$, there exists a prime $L^H$-ideal $\mf r_{H,c}$ satisfying the following constraints.
    \begin{itemize}
        \item $\mf r_{H,c}$ has ideal class $c$.
        \item $\mf r_{H,c}$ splits completely in $L$.
    \end{itemize}
\end{lemma}
\begin{proof}
    This is an application of the Chebotarev density theorem. For brevity, set $F\coloneqq L^H$, and let $M$ be the Hilbert class field of $F$. 
    
    The main claim is that $M\cap L=F$. Certainly $M\cap L$ contains $F$. However, $M/F$ is unramified, and $(\zeta-1)\cap F$ is totally ramified in $L$, so the claim follows.
    Thus,
    \[\op{Gal}(ML/F)\cong\op{Gal}\big(M/F\big)\times\op{Gal}\big(L/F\big)\cong\op{Cl}F\times H.\]
    To finish, choose $g\in\op{Gal}(M/F)$ corresponding to $c\in\op{Cl}F$ and then use the Chebotarev density theorem to find a prime $F$-ideal $\mf r$ such that $\op{Frob}_\mf r=(g,1)$. This prime works.
\end{proof}

\subsubsection{Choosing Cocycles}
We now choose our representative of $\alpha_2\in\widehat H^2(G,\op{Hom}_\ZZ(\ZZ[V_L],\AA_L^\times))$
For this, we must find explicit $2$-cocycles to represent the various $i_{v(u)}\big(u_{L_{v(u)}/K_u}\big)$s for $u\in V_K$. Some of these will be easy. For example, for the finite place $v=v((p))=\mf P$, we can set
\[c_p(\sigma^i,\sigma^j)=x^{-\floor{\frac{i+j}n}}\]
to represent $u_{L_\mf P/K_{(p)}}\in\widehat H^2(G,L_\mf P^\times)$, so we set $\widetilde c_p\coloneqq i_\mf P(c_p)$.

Additionally, for $v=v(\infty)$, we set $c_\infty\in Z^2(G_v,L_v/K_\infty)$ by $c_\infty(\tau^i,\tau^j)\coloneqq(-1)^{\floor{(i+j)/2}}$, where $\tau\coloneqq\sigma^{n/2}$. In order to make $i_v(c_\infty)$ compatible with $i_\mf P(c_p)$, we recall
\[\xi_{v,\langle e\rangle}\equiv i_\mf P(x)\cdot i_v(-1)\pmod{N_{\langle\tau\rangle}\AA_L^\times},\]
so we set $\widetilde c_\infty(\tau^i,\tau^j)=(\xi_{v,\langle e\rangle}/i_\mf P(x))^{\floor{(i+j)/2}}$ so that $[\widetilde c_\infty]=[i_v(c_\infty)]$ in $H^2(G_v,\AA_L^\times)$.

Lastly, we let $u=(q)$ denote a finite unramified place of $V_K$, and let $\mf Q$ be the finite place associated to $v(u)$. For brevity, set $H\coloneqq G_v$ and $F\coloneqq L^H$, and note $\op{Gal}(L/F)=H=\langle\sigma^{k_q}\rangle$, where $k_q\in\ZZ$ is chosen so that $x^{k_q}\equiv q\pmod{p^\nu}$. Now, let $n_q$ denote the order of $q\pmod{p^\nu}$; using our chosen uniformizer $\varpi_u\in\mf Q$, we can define the $2$-cocycle by
\[\left(\sigma^{k_qi},\sigma^{k_qj}\right)\mapsto\varpi_u^{\floor{\frac{i+j}{n_q}}},\]
for $0\le i,j<n_q$ to  $u_{L_v/K_u}\in\widehat H^2(H,L_v^\times)$. It will be helpful to be able to change between generators, so we pick up the following lemma.
\begin{lemma} \label{lem:cocyclegenchange}
    Let $G=\langle\sigma\rangle$ be a finite cyclic group of order $n$. Further, suppose $k\in\ZZ$ has $\gcd(k,n)=1$. Then define $\chi,\chi_k\in Z^2(G,\ZZ)$ by
    \[\chi\left(\sigma^i,\sigma^j\right)\coloneqq\floor{\frac{i+j}n}\qquad\text{and}\qquad\chi_k\left(\sigma^{ki},\sigma^{kj}\right)\coloneqq\floor{\frac{i+j}n},\]
    where $0\le i,j<n$. Then $[\chi]=k[\chi_k]$ in $H^2(G,\ZZ)$.
\end{lemma}
\begin{proof}
    It is well-known that
    \begin{equation}
        (-\cup[\chi_k])\colon\widehat H^0(G,\ZZ)\to\widehat H^2(G,\ZZ) \qquad\text{and}\qquad\left(-\cup\left[\sigma^k\right]\right)\colon\widehat H^2(G,\ZZ)\to\widehat H^0(G,\ZZ)\label{eq:cycliccupiso}
    \end{equation}
    are inverse isomorphisms. As such, we can compute
    \[\chi\cup\sigma^k=\sum_{g\in G}\chi\left(g,\sigma^k\right)=\sum_{\ell=0}^{n-1}\chi\left(\sigma^\ell,\sigma^k\right)=k,\]
    so $[\chi]=[k]\cup[\chi_k]=k[\chi_k]$.
\end{proof}
We now apply \Cref{lem:cocyclegenchange} to our context, using its notation. Define $d_q\coloneqq\gcd(k_q,n_q)$ so that $[\chi_{d_q}]=(k_q/d_q)[\chi_{k_q}]$. As such, we find $y_q\in\ZZ$ with $y_q\cdot k_q/d_q\equiv1\pmod{n_q}$ so that we can represent $\alpha(L_{v}/K_u)$ by the $2$-cocycle
\[c_q\colon\left(\sigma^{d_qi},\sigma^{d_qj}\right)\mapsto\varpi_u^{y_q\floor{\frac{i+j}n}}.\]
It remains to make $i_v(c_q)$ compatible with $i_\mf P(c_p)$. This is difficult, so we outsource to the following lemma.
\begin{lemma} \label{lem:constructgeneralxis}
    Fix everything as above. Then there exists $\xi_u\in L^\times$ and $e_u\in\ZZ$ such that
    \[\xi_u\varpi_u\equiv i_v(\varpi_u)\cdot i_\mf P(x)^{e_u}\pmod{N_H\AA_L^\times}.\]
\end{lemma}
\begin{proof}
    We
    begin by embedding $\varpi_u\in L^\times$ to $\AA_L^\times$, yielding
    \[\varpi_u=\prod_{w\in V_L}i_w(\varpi_u).\]
    We are interested in $\varpi_u\pmod{N_H\AA_L^\times}$. Most of these places $w\in V_L$ can be erased immediately: if $v'\in V_L$ is a finite place not lying over $(p)$, $\mf q$, nor $\mf r$, then we claim
    \[\prod_{w\in Hv'}i_w(\varpi_u)\]
    is a norm in $N_{H}\AA_L^\times$. Indeed, let $u'$ be the place of $V_F$ under $v'$. The extension $L_{v'}/F_{u'}$ is unramified because $v'$ does not lie over $(p)$, and $\varpi_u\in\OO_{F_{u'}}^\times$ because $v'$ does not lie over $\mf q$ or $\mf r$. Thus, there is some $a_{v'}\in L_{v'}$ such that $\varpi_u=\op N_{L_{v'}/F_{u'}}(a_{v'})=N_{H_{v'}}(a_{v'})$, so
    \[N_H(i_{v'}(a_{v'}))=\prod_{w\in Hv'}i_w(\varpi_u),\]
    follows by direct expansion. Collecting our norms, we see
    \begin{equation}
        \varpi_u\equiv i_v(\varpi_u)\cdot i_\mf P(\varpi_u)\cdot\prod_{w\mid\mf r}i_w(\varpi_u)\cdot\prod_{w\mid\infty}i_w(\varpi_u)\pmod{N_H\AA_L^\times}. \label{eq:varpiinitial}
    \end{equation}
    We now deal with each term individually.

    Fixing some place $v'\in V_L$ above $\mf r$, the fact that $\mf r$ is totally split in $L$ implies
    \[\prod_{w\mid\mf r}i_w(\varpi_u)=\prod_{h\in H}i_{hv'}(\varpi_u)=N_H(i_{v'}(\varpi_u))\in N_H\AA_L^\times,\]
    so we may ignore these terms. Next up, we deal with the finite place $\mf P$ above $(p)$. Set $\mf p\coloneqq F\cap\mf P$, and note $\OO_F/\mf p\cong\FF_p\cong\OO_K/(p)$, and this isomorphism is unique and so induced by $K\subseteq F$, so $\varpi_u\pmod{\mf p}$ is a power of $x\pmod{\mf p}$. Namely, there is some power $e'$ such that $\varpi_u/x^{e'}\equiv1\pmod{\mf p}$, which is a norm from $L_\mf P$. Thus,
    \[i_\mf P(\varpi_u)\equiv i_\mf P(x)^{e'}\pmod{N_H\AA_L^\times}.\]
    It remains to deal with infinite places, which is where $\xi_u$ will come from. We begin by fixing some infinite place $v'\mid\infty$. We have two cases.
    \begin{itemize}
        \item If $\tau\notin H$, then we see that
        \[\prod_{w\mid\infty}i_w(\varpi_u)=\prod_{h\in H}i_{hv'}(h\varpi_u)=N_H(i_{v'}(\varpi_u))\in N_H\AA_L^\times,\]
        where the last step is because $H_{v'}=\{1\}$. So we may set $\xi_u=1$ and $e_u=0$ to finish. 
        \item Otherwise, $\tau\in H$. For concreteness, associate $v'$ to the embedding $\sigma\colon L\to\CC$. In fact, $\sigma\left(L^H\right)\subseteq\RR$ because $L^H$ is fixed by $\tau\in H$, so we'll consider
        \[i_{v'}\left(\sqrt{\varepsilon_{u,v'}\sigma(\varpi_u)}\right)\in\AA_L^\times,\]
        where the sign $\varepsilon_{u,v'}\in\{\pm1\}$ is chosen to ensure $\varepsilon_{u,v'}\sigma(\varpi_u)>0$. Concretely, $\sqrt{\varepsilon_{u,v'}\sigma(\varpi_u)}$ is a Cauchy sequence of elements of $L^H$ under the metric induced by $\sigma\colon L^H\to\RR$, whose square approaches $\varepsilon_{u,v'}\sigma(\varpi_u)>0$. Applying any $h\colon L_{v'}\to L_{hv'}$ for $h\in H$ will not change that this sequence is Cauchy, nor will it change the limit, because all relevant elements are in $L^H$. Thus, without ambiguity, we may compute
        \begin{align*}
            N_Hi_{v'}\left(\sqrt{\sigma(\varepsilon_{u,v'}\varpi_u)}\right) 
            &= \prod_{h\langle\tau\rangle\in H/\langle\tau\rangle}i_{hv'}\left(\sqrt{\sigma(\varepsilon_{u,v'}\varpi_u)}\cdot\tau\sqrt{\sigma(\varepsilon_{u,v'}\varpi_u)}\right)
            = \prod_{w\in Hv'}i_w(\varepsilon_{u,v'}\varpi_u).
        \end{align*}
        In total, we see that
        \[\prod_{w\in Hv'}i_w(\varpi_u)\equiv\prod_{w\in Hv'}i_w(\varepsilon_{u,v'})\equiv\left(\xi_{v',H}\cdot i_\mf P(x)^{-|H|/2}\right)^{(1-\varepsilon_{u,v'})/2}\pmod{N_H\AA_L^\times}\]
        by \Cref{lem:constructgeneralxis}. Now, $H$ acts on the places above $\infty$, so choosing representatives $v'$ from the orbits of this action and collecting the above products completes the proof.
        \qedhere
    \end{itemize}
\end{proof}
Lastly, we fix the $i_\mf P$ term. For this, we use the following lemma.
\begin{lemma} \label{lem:fixpplace}
    Fix everything as above. Suppose that we have a subgroup $H\subseteq G$ and power $e\in\ZZ$ such that
    \[[i_\mf P(x)^e]=[1]\]
    as elements of $\widehat H^0(H,\AA_L^\times/L^\times)$. Then $i_\mf P(x)^e\equiv1\pmod{N_H\AA_L^\times}$.
\end{lemma}
\begin{proof}
    The point is to show that $|H|$ divides $e$. Let $H=\langle\sigma^d\rangle$ for a fixed $d\mid n$. We note
    \[\op{Res}\widetilde c_p\colon\left(\sigma^{di},\sigma^{dj}\right)\mapsto i_\mf P(x)^{-\floor{\frac{i+j}{n/d}}}\]
    represents the fundamental class of $\widehat H^2(H,\AA_L^\times/L^\times)\simeq\ZZ/|H|\ZZ$, so $i_\mf P(x)^{-1}=(\op{Res}\widetilde c_p)\cup\sigma^d$ is a generator of $\widehat H^0(H,\AA_L^\times/L^\times)$ and so has order $|H|$.
    Thus, $[i_\mf P(x)]^e=[1]$ implies that $|H|$ divides $e$. It follows $N_H(i_\mf P(x))^{e/|H|}=i_\mf P(x)^e$ is a norm in $N_H\AA_L^\times$, which is what we wanted.
\end{proof}
We now return to our computation. Currently, we have some $\xi_u$ and $e_u$ such that
\[\xi_u\varpi_u\equiv i_v(\varpi_u)\cdot i_\mf P(x)^{e_u}\pmod{N_H\AA_L^\times}.\]
However, we know abstractly that the $2$-cocycles $i_v(c_q)$ and $\op{Res}\widetilde c_p$ both represent the fundamental class of $\widehat H^2(H,\AA_L^\times/L^\times)$, so
\[\left[i_v(\varpi_u)^{y_q}\right]=[i_v(c_q)]\cup\left[\sigma^{d_q}\right]=[\op{Res}\widetilde c_p]\cup\left[\sigma^{d_q}\right]=\left[i_\mf P(x)^{-1}\right]\]
as elements of $\widehat H^0(H,\AA_L^\times/L^\times)$. Combining,
\[[1]=[\xi_u\varpi_u]^{y_q}=\left[i_v(\varpi_u)^{y_q}\cdot i_\mf P(x)^{y_qe_u}\right]=\left[i_\mf P(x)^{y_qe_u-1}\right]=[i_\mf P(x)]^{y_qe_u-1}\]
as elements of $\widehat H^0(H,\AA_L^\times/L^\times)$. Thus, \Cref{lem:fixpplace} lets us conclude that $i_\mf Px^{y_qe_u}\equiv i_\mf Px\pmod{N_H\AA_L^\times}$, so
\[(\xi_u\varpi_u)^{y_q}\equiv i_v(\varpi_u)^{y_q}\cdot i_\mf P(x)\pmod{N_H\AA_L^\times}.\]
In total, we can choose
\[\widetilde c_q\left(\sigma^{di},\sigma^{dj}\right)\coloneqq\left(\xi_u^{y_q}\varpi_u^{y_q}/i_\mf P(x)\right)^{\floor{\frac{i+j}{n_q}}}\]
to represent $i_v(u_{L_v/K_u})\in\widehat H^2(H,\AA_L^\times)$.

To synthesize all places, we set
\begin{equation}
    \omega_u\coloneqq\begin{cases}
        1 & \text{if }u=(p), \\
        \xi_\infty & \text{if }u=\infty, \\
        \xi_u^{y_q}\varpi_u^{y_q} & \text{if }u\notin\{(p),\infty\},
    \end{cases}\qquad\text{and}\qquad \text{if }d_u\coloneqq\begin{cases}
        1 & \text{if }u=(p), \\
        n/2 & \text{if }u=\infty, \\
        d_q & \text{if }u\notin\{(p),\infty\},
    \end{cases} \label{eq:omegadef}
\end{equation}
so that
\[\widetilde c_u\left(\sigma^{d_ui},\sigma^{d_uj}\right)=(\omega_u/i_\mf P(x))^{\floor{\frac{d_ui+d_uj}{n}}}\]
in all cases.

\subsection{Finishing the Computation} \label{subsec:compute-a3}
In this subsection we continue in the context of \cref{subsec:compute-a2} and complete the computation of the Tate canonical class. To begin, set $\alpha^{(1)}\in\left(\AA_L^\times/L^\times\right)^G$ by $\alpha^{(1)}\coloneqq i_\mf P(x)^{-1}$ so that the corresponding $2$-cocycle $c^{(1)}\in H^2(G,\AA_L^\times/L^\times)$ represents the global fundamental class.

To compute the class $\alpha_2\in H^2(G,\op{Hom}_\ZZ(\ZZ[V_L],\AA_L^\times))$, it suffices by \Cref{prop:tuple-to-cocycle}, it suffices to define the element $\alpha^{(2)}\in\op{Hom}_{\ZZ[G]}(\ZZ[V_L],\AA_L^\times)$ by
\[\alpha^{(2)}\colon gv(u)\mapsto g\omega_u/i_\mf P(x)\]
for each $u\in V_K$ and $g\in G$. Indeed, let $c^{(2)}$ be the corresponding $2$-cocycle. Tracking through the definition of $\alpha_2$ via \Cref{lem:magicaltate}, we see that
\[\left(\op{eval}_{v(u)}\op{Res}_{G_{v(u)}}c_2\right)\left(\sigma^{d_ui},\sigma^{d_uj}\right)=(\omega_u/i_\mf P(x))^{\floor{\frac{d_ui+d_uj}n}}=\widetilde c_u\left(\sigma^{d_ui},\sigma^{d_uj}\right),\]
so $c^{(2)}$ does indeed represent $\alpha_2$.

Lastly, we define $\alpha^{(3)}\in\op{Hom}_{\ZZ[G]}(\ZZ[V_L]_0,L^\times)$ by
\[\alpha^{(3)}\colon(gv(u)-\mf P)\mapsto g\omega_u,\]
and let $c^{(3)}\in H^2(G,\op{Hom}_\ZZ(\ZZ[V_L]_0,L^\times))$ be the corresponding $2$-cocycle. A computation shows that the diagram
\[\begin{tikzcd}
	0 & {\ZZ[V_L]_0} & {\ZZ[V_L]} & \ZZ & 0 \\
	0 & {L^\times} & {\AA_L^\times} & {\AA_L^\times/L^\times} & 0
	\arrow[from=1-1, to=1-2]
	\arrow[from=1-2, to=1-3]
	\arrow[from=1-3, to=1-4]
	\arrow[from=1-4, to=1-5]
	\arrow[from=2-1, to=2-2]
	\arrow[from=2-2, to=2-3]
	\arrow[from=2-3, to=2-4]
	\arrow[from=2-4, to=2-5]
	\arrow["{\alpha^{(2)}}", from=1-3, to=2-3]
	\arrow["{\alpha^{(1)}}", from=1-4, to=2-4]
	\arrow["{\alpha^{(3)}}", from=1-2, to=2-2]
\end{tikzcd}\]
commutes, so $c^{(3)}$ indeed represents the canonical class $\alpha_3$. This is what we wanted.

\printbibliography[title={References}]

\end{document}